\documentclass{svjour33}
\usepackage{}
\usepackage{pifont}                     
\smartqed  

\usepackage{graphicx}
\usepackage{amsmath}
\usepackage{amssymb}
\usepackage{latexsym}
\usepackage{subfigure}

\newcommand{\be}{\begin{equation}} \newcommand{\ee}{\end{equation}}
\newcommand{\barr}{\begin{array}} \newcommand{\earr}{\end{array}}
\newcommand{\bitem}{\begin{itemize}}
\newcommand{\eitem}{\end{itemize}}

\newtheorem{thm}{Theorem}[section]
\newtheorem{alg}{Algorithm} 

\begin{document}

 \mag=1200
\title{ Orthogonally Accumulated Projection Methods for Linear System of Equations}

\author{ Wujian Peng \and  Shuhua Zhang}

\institute{Wujian Peng \at
              Dept. of Math. and stat. Zhaoqing Univ. China, 526061
              \email{wpeng@zqu.edu.cn}           
           \and
            Shuhua Zhang \at
            Dept. of Math.,  Tianjin Univ. of Finance and Economics
           \email{szhang@tjufe.edu.cn}
}

\date{Received: date / Accepted: date}

\maketitle
\begin{abstract}
    A type of iterative orthogonally accumulated projection methods for solving linear system of equations are proposed in this paper. This type of methods are applications of accumulated projection(AP) technique proposed recently by authors. Instead of searching projections in a sequence of subspaces as done in the original AP approach, these methods try to efficiently construct a sequence of orthonormal vectors while the inner-product between the solution to the system and each vector in the sequence can be easily calculated, thus the solution can be retrieved in finite number of iterations in case of exact arithmetic operations. We also discuss the strategies to handle loss-of-orthogonality during the process of constructing orthonormal vectors. Numerical experiments are provided to demonstrate the efficiency of these methods.

\end{abstract}

\section{Introduction}
  The study of iterative methods for solving least square problems in the form
  \be \label{def:sys}
     Ax = b
   \ee
   where $A\in R^{n\times n}$, especially for large scale computing is of vital importance. Here we always assume   $A$ is nonsingular so that there exists a unique solution to the system. There are a lot of iterative methods available\cite{Axelsson1,Axelsson2,Templates,Golub,Hackbusch_MG,Hackbusch_It} for solving system (\ref{def:sys}).

  Recently all current iterative methods are classified as extended Krylov subspace methods in \cite{pengPAP}, which are characterized by their major operations: matrix-vector multiplications with usually  one or two fixed matrices and one or two fixed initial vectors. These includes the most well-known stationary methods such as Jacobi, Gauss Seidal as well as SOR methods with their iterative matrices formed on the base of splitting the coefficient matrices\cite{Golub}, and the row projection methods such as Karcmarz's method and Cimmino's methods where the iterative matrix (not explicitly formed in iterations) are constructed by the successive multiplications of a sequence of projection matrices with a fixed sequence length $m$(depending on the splitting of the coefficient matrix into $m$ submatrices\cite{Bramley}\cite{Galantai}). The non-stationary iterative methods include the well-known Krylov subspace methods such as conjugate gradient method(CG) for symmetric positive definite systems, MINRES, SYMMLQ for general symmetric but indefinite systems, and GMRES, BiCG, BiCR, QMR, LSQR, etc. for general nonsymmetric systems\cite{Freund1991,Golub,Paige,Saad,Vorst}; many of these methods(including GMRES, MINRES, SYMMLQ, MINRES, QMR, LSQR) use the strategy of reducing some related residual norms to search for approximate solutions, while variants of CG and BiCG methods use the strategy of producing a sequence of orthogonal residuals, thus they can reach the exact solutions with $n$ iterations in exact arithmetic operations, where $n$ is the number of unknowns\cite{Golub}.

  In \cite{pengPAP} authors also presents a first non-krylov subspace type methods--The Accumulated Projection Methods. These type of methods rely on successive projections over subspaces of $R^n$, which produce a sequence of projection vectors with a monotonically increasing Euclidean norms. Unlike the well-known row-projection technique which can be shown as a traditional stationary iterative methods \cite{Galantai}, the AP methods proposed in\cite{pengPAP} do not involve matrix-vector multiplications with any fixed matrices and fixed vectors. Equipped with some accelerating technique, the AP methods exhibit some superior behavior than traditional extended Krylov subspace methods\cite{pengPAP} in some cases.

  The success of AP methods rely on the calculation of projection vector of exact solution $x\in R^n$ over a sequence of subspaces $W_k$ of $R^n$ $(k=1,2, \cdots, m)$, where $W_k$ is formed by the row vectors of coefficient matrix and the most recent approximations $p_k$ of $x$. The calculation of these projection vectors are based on the QR factorization of matrix $W_k$ for general matrices, or QS\cite{pengLGO} decomposition of $W_k$ if the coefficient matrix is sparse. Generally speaking, the QR factorization needs $O(m^2n)$ flops and is thus a heavy burden if a long iteration is needed, current LGO decomposition requires that the coefficient matrix satisfies some special property (for example, $k$-orthogonality) and its implementation is quite complicate. One of our purpose in this paper is to provide a more efficient way to handle the projection of any given vector into a subspace of $R^n$ with much less float point operations.

  Our major task in this paper is to provide a class of methods based on the principle of accumulated projection to handle linear system of equations. For the sake of completeness, we are to briefly review the principle of  accumulated projection technique and its applications in the next section. The other sections are devoted to discuss the exploration of AP technique in a more intricate way which leads to a series of algorithms for solving linear systems.

\section{Principle of AP technique}

Now we review the basic idea of accumulated projection methods.  To approximate any vector $x$ in $R^n$, one has to construct a subspace $W$ of $R^n$ with a much smaller rank than $n$  so that a ``projection" vector $p$ of $x$ is easily available.  Current prevalent methods depend on the strategy of reducing the length of residual vectors to obtain such a projection. While only a few methods use the regular orthogonal projection to get approximate vectors, which include the so-called General Error Minimizing Method (which is similar to GMRES method)\cite{Rainald} and the Line Projection method proposed in \cite{pengLP}, both can be classified as extended Krylov subspace methods since both of them depend on certain Krylov subspace from which a projection vector is sought. To be able to figure out the projections of $x$ over subspace $W$, one has to get some ``footprint" of $x$ over $W$,  for example in GMRES-like methods a basis vectors of $W$ in the form of $A^k b$ with $b$ as image of $x$ under the transformation $A$ are required, while in GMERR and LP methods, \textbf{the inner-products between $x$ and a basis of $W$ are available.} By this observation we can derive another class of methods for solving linear system of equations using orthogonal projections.

The basic idea of AP is to use the orthogonal projections of vector $x$ as its approximations, while each projection is used to form another subspace from which a better approximation is sought. The following graph can be used to illustrate the whole idea.
\begin{figure}[h]
\begin{center}
\includegraphics[height=6cm,width=8cm]{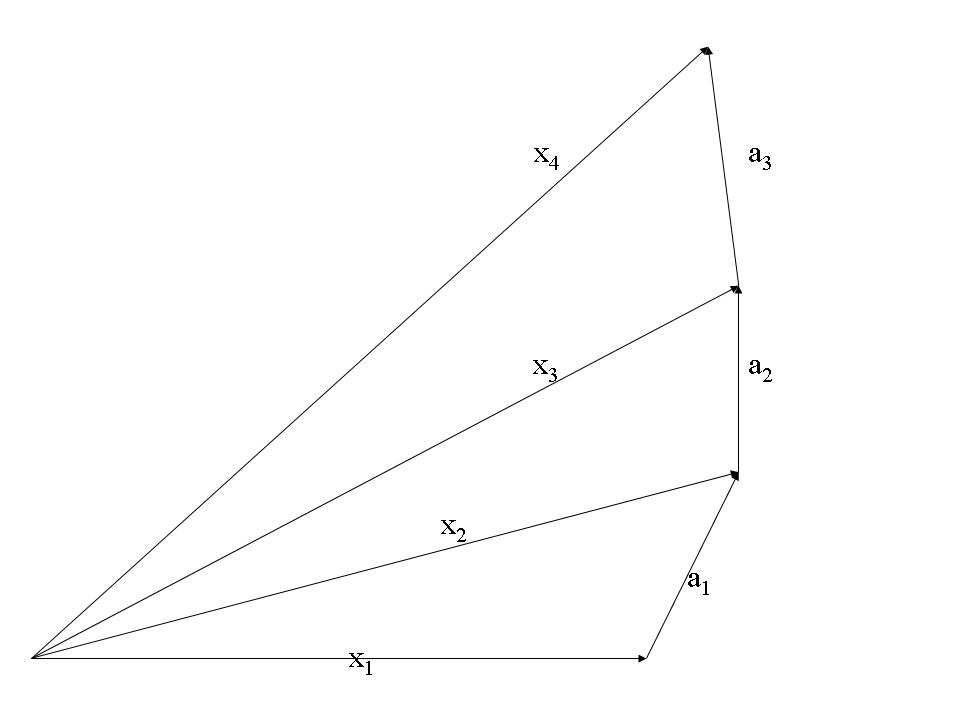}
\caption{Accumulated Projection}
\label{fig:AP_principle}
\end{center}
\end{figure}
where $x_i$ stands for the approximations to $x$ and $a_i$ are projection vectors of $x$ on some subspaces of $R^n$. $x_{i+1}$ is the projection of vector $x$ in a subspace $W_i$ formed by $x_i$ and a subspace $\tilde{W}_i$ where projection vector $a_i$ of vector $x$ is easily available.

The following algorithm describes a simple implementation of the accumulated projection idea, where vector $a_i$ is orthogonal to vector $x_i$.

\begin{alg}\label{alg:AP}(accumulated projection process-AP) The following procedure
produces an approximate vector $p$ to the solution vector $x$ which satisfies $Ax = b.$
\begin{quote}
\begin{itemize}
\item[(1)]Divide matrix $A$ into $k$ blocks:
 $A  = [A_1'. A_2', \cdots, A_k']'$, divide $b$  correspondingly:  $b = (b_1', b_2', \cdots,
 b_k')'$.
 \item[(2)] Initialize $p_0$ as $ p_0= \alpha A' b$ and $c_0 = \alpha ||b||^2$, where \mbox{$\alpha = ||b||^2 / ||A'b||^2$}.
\item[(3)] For $i =1$ to $k$
 \begin{quote}
\begin{itemize}
\item[(3.1)] Construct matrix $W_i =[p_{i-1}, A_i'] $ and vector $l =[c_{i-1},b_i']'$.
\item[(3.2)] Compute the projection vector $p_i$ of $x$ onto subspace $ran(W_i)$  and the scalar $c_i(=x'p_i) $.
\end{itemize}
\end{quote}
\item[(4)] Output $p(=p_k)$ and $c(=c_k)$.
\end{itemize}
\end{quote}
\end{alg}

This algorithm formed the basis of some more efficient solvers for linear system of equations such as SAP and MSAP and APAP methods introduced in \cite{pengSAP} and \cite{pengPAP}. It is observed that these methods seem to be more efficient than regular Krylov subspace methods in case of large scale systems in some situation. It is necessary to mention that these methods do not construct any Krylov subspace methods and thus can not be classified as extended Krylov subspace methods.  In this paper we will show that the AP process can also be used to construct a class of Krylov subspace methods, named as orthogonally accumulated projection solver(OAP).

\section{An orthogonally accumulated projection through tridiagonalization}\label{sec:tridiag}

In this section we will consider to solve system (\ref{def:sys}) with a unsymmetric coefficient matrix $A$. The main idea is to transform the original system (\ref{def:sys}) into a system
     \be
          Qx = c
     \ee
where $Q$ is an orthogonal matrix, i.e., $Q'Q = I$ where $I$ is the identity matrix. In other words, we will search for a sequence of orthonormal  vectors $v_i \,(i=1,2,\cdots, n)$ and real numbers $c_i \, (i=1,2,\cdots, n)$ so that $x'v_i = c_i$, and thus $x $ can be taken  as  $\sum_i^n c_i v_i$. In the meantime we do not have to spend too much extra storage space to store all vectors $v_i(i=1,2,\cdots, n)$, instead we will show that a short length recurrence relationship occurs between contagious orthogonal vectors so that only a few extra storage space for these vectors is needed.

In order to figure out how this will work, let us review the principle of AP as illustrated in Figure (\ref{fig:AP_principle}). In general the sequence of projection vectors $a_i$ come from some predetermined subspaces and thus they are not necessary to be orthogonal. However it is possible for us to work out a way so that all of these projection vectors $a_i \,(i=1,2, \cdots, )$ form an orthogonal sequence. To be complete, we first recall the Laczos iterations for tridiagonalization of a rectangular matrix.

\subsection{Matrix tridiagonalization by orthogonal transformation}

Any matrix $A \in R^{n\times n}$ can be transformed into the following tridiagonal form
  \begin{equation} \label{eq:tridiagonal_form1}
     U'AV = T
  \end{equation}
where $T$ is tridiagonal
 $$ T= \left(\begin{array}{cccccc}
            \alpha_1 & \beta_1  & 0       & \cdots&  0 & 0 \\
            \gamma_1 & \alpha_2 & \beta_2 & \cdots& 0 & 0 \\
            \vdots   & \vdots &  \vdots &   \vdots & \vdots & \vdots  \\
            0        &   0    &  0      &    \cdots & \alpha_{n-1}& \beta_{n-1} \\
             0        &   0    &  0      &    \cdots & \gamma_{n-1}& \alpha_n \\
             \end{array} \right),
$$ both $U$ and $V$ are orthogonal, i.e., $U'U  =V'V = I_m$. This transform can be accomplished in a rather stable way by applying Householder transformations on both sides of $A$. However when $A$ is sparse and large, we can expect dense and large submatrices to appear in this process, which makes it not suitable in large scale computations.

Fortunately a Lanczos-like process can be used to do the tridiagonalization in a much cheaper and efficient way. To illustrate this we rewrite equation (\ref{eq:tridiagonal_form1}) into the following forms
 \begin{equation} \label{eq:tridiagonal_form2_1}
    A V = U T
   \end{equation}
and
 \begin{equation} \label{eq:tridiagonal_form2_2}
    A' U= V T'
   \end{equation}
Equating $k$-th column of both sides of (\ref{eq:tridiagonal_form2_1}) and (\ref{eq:tridiagonal_form2_2})
we have
 \begin{equation} \label{eq:tridiagonal_form2_3}
    A v_k = \gamma_{k} u_{k+1} + \alpha_k u_k + \beta_{k-1} v_{k-1}
   \end{equation}
and
 \begin{equation} \label{eq:tridiagonal_form2_4}
    A' u_k = \beta_k v_{k+1} + \alpha_k v_k + \gamma_{k-1} v_{k-1}
   \end{equation}
with $ \beta_0 = \gamma_0 = 0$, where $v_k$ and $u_k$ denote the $k$-th columns of
matrix $V$ and $U$ separately, $ v_0$ and $ u_0 = 0$ are zero vectors, i.e, $V = [v_1, v_2, \cdots, v_n]$ and
$U = [u_1, u_2, \cdots , u_n]$. Especially we have
\begin{equation} \label{eq:1st_col}
   A v_1 = \alpha_1 u_1 + \gamma_1 u_2 \mbox{~~~ and~~~~ } A' u_1 = \alpha_1 v_1 + \beta_1 v_2
\end{equation}
which suggests that if both $u_1$ and $v_1$ are given, then $v_2$ and $u_2$ can be calculated simultaneously.
The rest vectors $v_k$ and $u_k$ for $k\ge 3$ can be calculated by rewriting (\ref{eq:tridiagonal_form2_3})
and (\ref{eq:tridiagonal_form2_4}) as follows
\begin{equation} \label{eq:tridiagonal_form2_5}
    u_{k+1} = \frac{1}{\gamma_{k}} (A v_k - \alpha_k u_k - \beta_{k-1} v_{k-1} )
   \end{equation}
and
 \begin{equation} \label{eq:tridiagonal_form2_6}
   v_{k+1} = \frac{1}{\beta_k} (A' u_k - \alpha_k v_k - \gamma_{k-1} v_{k-1})
   \end{equation}
The following algorithm depicts the above process.

\begin{alg}\label{alg:OAP0}
   Let $A \in R^{n \times n}$, $v_1 $ and $u_1 $ be  unit vectors.
  \begin{quote}
     $ \beta_0 = \gamma_0 = 0$, $v_0 = u_0 =0(\in R^n)$ \\
      for $k=1$ to $n-1$ \\
       \mbox{~~~}$\alpha_k = u_kAv_k$ \\
       \mbox{~~~}$w_k = Av_k - \alpha_k u_k -\beta_{k-1} u_{k-1}$\\
       \mbox{~~~}$\gamma_k = || w_k||$\\
       \mbox{~~~}$u_{k+1} = w_k / \gamma_k$\\
       \mbox{~~~}$q_k = A' u_k - \alpha_k v_k -\gamma_{k-1} v_{k-1}$\\
       \mbox{~~~}$\beta_k = ||q_k||$\\
       \mbox{~~~}$ v_{k+1} = q_k/ \beta_k$\\
      end
  \end{quote}
\end{alg}

In exact arithmetic operations the above Lanczos-like iteration will produce two orthonormal vector squences $v_1, v_2, \cdots v_n$ and $u_1, u_2, \cdots u_n$ with any starting unit vectors $u_1$ and $v_1$, assuming no break-down happens(i.e.,$ \gamma_k \neq 0$ and $ \beta_k \neq 0$ for all $k$). Note that in each loop in the iteration one needs only two matrix-vector multiplications as its major flop counts, this makes it very effective when dealing with tridiagonalizations of large and sparse matrices.

\subsection{orthogonally accumulated projection}

We have observed that in basic AP algorithm to make sure next approximation $x_{k+1}$ is a better approximation to $x$(the exact solution) than $x_k$, a projection on a subspace which contains $x_k$ must be done, which guarantees that
$ ||e_{k+1}||<||e_k || $ where $e_k = x-x_k$ is the error vector associated with $x_k$. However if $a_i$ can be constructed in such a way that they always satisfy
\begin{equation} \label{eq:11}
x_k \perp a_i \mbox{ for } i>k, \end{equation}
there is no need to do the extra projection to get the next approximation $x_{k+1}$, instead one can simply obtain $x_{k+1}$ by $x_{k+1} = x_k + a_i$. Obviously if vector sequence $\{a_k\}_1^n$ forms an orthonormal sequence of vectors in $R^n$, and let $x_k  = \sum_{i=1}^k c_i a_i$ where  $c_i = x' a_i$, then it is easy to see that
(\ref{eq:11}) holds true. This is exactly the principle of orthogonally accumulated projection(OAP). In other words, to solve system (\ref{def:sys}), OAP method builds a sequence of orthonormal vectors $\{v_i\}_1^n$ as well as sequence of $\{c_i\}_1^n$, the inner-product between $x$ and each of $v_i$, i.e., $c_i = x'v_i$, thus $x$ can be retrieved as  $x = \sum_{i=1}^n c_i v_i$.

We will shown in next section that in exact arithmetic operations,  Algorithm \ref{alg:OAP0} will produce a sequence of orthonormal vectors $\{v_i\}_{i=1}^n$; in order to find the inner-product $c_i$ between $x$ and each $v_i$, we multiply by $x$ both sides of equation (\ref{eq:tridiagonal_form2_6}), this leads
to
\be\label{comput:ci}
\begin{array}{ll}
    c_{k+1} &= \frac{1}{\beta_k} ( x'A'u_k  - \alpha_k x'v_k - \gamma_{k-1}x'v_{k-1}) \\[3mm]
                       &=  \frac{1}{\beta_k} (b'u_k - \alpha_k c_k - \gamma_{k-1} c_{k-1})
\end{array}
\ee
since $Ax = b$, particularly we have $c_2 = \frac{1}{\beta_1} (b'u_1 - \alpha_1 c_1 ) $.
This implies that if $c_1$ is known, then all the other subsequent $c_i(i=2,3,\cdots n)$ can be calculated by (\ref{comput:ci}).
These process can be described in the following algorithm, which is called orthogonally accumulated projection for solving linear system of equations.

\begin{alg}\label{alg:OAP3}(orthogonally accumulated projection method-OAP) Let $A$ in $R^{n\times n}$ be an unsymmetric and nonsingular matrix and $b\in R^n$ a non-zero vector. Let $v_1$ and $u_1$ be two unit vectors and $c_1 = x'v_1$ be given, where $x$ is the solution to (\ref{def:sys}).
The following process gives the exact solution $x$ to system $Ax = b.$
   \begin{quote}
      $ x_1 = c_1 v_1, \beta_0 = \gamma_0 = 0$, $v_0 = u_0 =0(\in R^n)$\\
      for $k =1 $ to $n - 1$\\[-4mm]
        \begin{quote}
           $\alpha_k = u_k' A v_k$ \\
           $ p_k = Av_k - \alpha_k u_k - \beta_{k-1} u_{k-1} $\\
           $\gamma_k = ||p_k||$ \\
           $u_{k+1} = p_k / \gamma_k $ \\
           $q_k = A'u_k - \alpha_k v_k - \gamma_{k-1} v_{k-1} $\\
           $\beta_k = || q_k || $\\
           $ v_{k+1} = q_k / \beta_k $\\
           $c_{k+1} = \frac{1}{\beta_k}(b'u_k - \alpha c_k - \gamma_{k-1} c_{k-1}) $\\
           $x_{k+1} = x_k + c_{k+1} v_{k+1}$ \\
        \end{quote}\vskip -4mm
      end
   \end{quote}
\end{alg}
Note that there are only two matrix-vector multiplications involved, and storage for extra four vectors
is needed besides that for the coefficient matrix $A$. In case $A$ is sparse(having an average of $m$ none-zero elements in each row) and large, the flop counts for one sweep of the loop is $O(mn)$. Therefore in exact arithmetic operations, there are only $O(m^2 n^2)$ flops needed for the whole procedure.

Remark: there are many options for the initial vectors $v_1$ while $u_1$ can be chosen arbitrarily. For example any row vector $A_i$ of matrix $A$ can be used for constructing $v_1$  ($v_1 = A_i'/ ||A_i|| $ with $c_1 = b_1/||A_i||$. Another type of options is any vector in the form $v_1 = t A'w$ where $w$ is any none-zero vector and $t$ is a scalar such that $v_1$ is a unit vector, and in this case one can see that $c_1$ can be obtained as
$ c_1 = t b'w$.

\subsection{Analysis of OAP}

In this section we discuss some properties of OAP as a direct method(in exact arithmetic operations). Note that
any unsymmetric matrix can also be transformed by Householder transformation into tridiagonal matrix $T$($T= V'AV$) with $V$ as orthogonal matrix, which suggests us to develop a similar algorithm for this type of transformation. However it turns out such a Lanczos-like iteration does not exist at least for arbitrarily chosen initial unit vector $v_1$. It is thus necessary to verify the orthonormality of the vectors sequences $\{v_i\}_1^n$ and $\{u_i\}_1^n$ in Algorithm \ref{alg:OAP3}.

\begin{thm}
Let $A$ be unsymmetric and nonsingular, $b\in R^n$ and $x$ is the solution to $Ax=b$. The vector sequence $v_k \, (k=1,2,\cdots, n)$ and $u_k\, (k=1,2,\cdots, n)$produced in Algorithm \ref{alg:OAP3} are orthonormal, assuming no breakdown happens, i.e., $\beta_k \neq 0$ and $\gamma_k \neq 0$ for any $k=1,2,3,\cdots, n-1$.
\end{thm}

Proof. Apparently all vectors $v_i$ and $u_i$ ($i = 1, 2, \cdots, n$) are unit vectors.
We first  show that $v_2'v_1=0$ and $u_2 'u_1=0$.

Note that
\begin{eqnarray}
     v_2'v_1 = 0  & \Leftrightarrow (A'u_1 - \alpha_1 v_1)' v_1 = 0
                  &\Leftrightarrow \alpha_1 = v_1'A'u_1
\end{eqnarray}
the last equation is exactly how $\alpha_1$ is constructed in the algorithm, hence we have $v_2$ and $v_1$ are orthogonal. Similary we have
\begin{eqnarray}
     u_2'u_1 = 0  & \Leftrightarrow (Av_1 - \alpha_1 u_1)' u_1 = 0
                  &\Leftrightarrow \alpha_1 = v_1'A'u_1
\end{eqnarray}
By induction, we assume $v_1, v_2, \cdots, v_k$ and $u_1, u_2, \cdots, u_k$ are orthonormal sequences of vectors, we need to show that $v_{k+1}'v_i =0$ and $u_{k+1}'u_i =0$ for $i \le k$.

In fact $$ u_{k+1}' u_k = 0 \Leftrightarrow (Av_k - \alpha_k u_k - \beta_{k-1} u_{k-1} )' u_k =0
           \Leftrightarrow \alpha_k = u_k' A v_k$$
        and
$$
\begin{array}{ll}
           u_{k+1}' u_{k-1} =0 & \Leftrightarrow (A v_k - \alpha_k u_k - \beta_{k-1} u_{k-1})' u_{k-1} =0\\
                &\Leftrightarrow \beta_{k-1} = u_{k-1}' A v_k \\
               &\Leftrightarrow \beta_{k-1} = (A' u_{k-1})' v_k \\
               &\Leftrightarrow \beta_{k-1} = v_k' (A' u_{k-1})\\
              & \Leftrightarrow \beta_{k-1} = v_k' ( \beta_{k-1} v_k - \alpha_{k-1} v_{k-1} - \gamma_{k-2} v_{k-2})\\
                &\Leftrightarrow \beta_{k-1} = \beta_{k-1}
\end{array}
$$
For $i \le k-2 $ we have
$$
\begin{array}{ll}
    u_{k+1}' u_i =0 &\Leftrightarrow (A v_k - \alpha_k u_k - \beta_{k-1} u_{k-1})' u_i =0  \\
                   &\Leftrightarrow  v_k' (A'u_i) =0 \\
                   &\Leftrightarrow  v_k' ( \beta_i v_{i+1} + \alpha_i v_i + \gamma_{i-1} v_{i-1} =0
\end{array}
$$ The last equation holds true since by assumption we have $v_k$ are orthogonal to $v_i$ for any $i<k$.
Similarly one can prove $u_{k+1}'u_i=0$ for $i \le k$.  $\Box$

\subsection{Control of loss of orthogonality}

There are several well-known Krylov subspace methods based on Lanczos iterations. The most famous method might the the wide-spread conjugate gradient method(CG)(by  Hestenes and Stiefel). Other effective methods include MINRES, SYMMLQ and LSQR(by Paige and Saunders), BiCG(by Fletcher) and BiCGstab(by Van der Vorst) and QMR(by Freund and Nachtigal), etc.  All of these methods(except CG) adopt the strategy of minimizing certain type of residual norm in related Krylov subspace.

Unfortunately Lanczos process often suffers severe loss of orthogonality, which explains the possible instability of most of the above Krylov subspace methods based on Lanczos iteration. It seems that there is no effective way to handle this issue in general. Krylov subspace methods based on Arnoldi iteration(such as GMRES) seems to be more stable but they usually need more storage requirement and flops in each iteration and thus usually have to be restarted.

Krylov subspace methods based on minimizing residual norms usually ignore the issue of loss of orthogonality. However it is vital to our orthogonally accumulated projection method. Fortunately we have an easy approach to detect whenever loss of orthogonality happens. Our approach is to make sure in every iteration the ``accumulated" vector $a_{k+1}$ is guaranteed to be orthogonal to current approximation $x_k$. Note that $x_k$ is a linear combination of $v_1, v_2, \cdots, v_k$ and $a_{k+1}= c_{i+1}v_{k+1}$ (with $c_{i+1}$ a real number) is supposed to be orthogonal to all $v_i$ for $i\le k$. Thus the angle between $x_k$ and $v_{k+1}$ a is good indicator when loss of orthogonality occurs. And whenever loss of orthogonality happens, we restart the OAP process on the residual equation $r_k = A e_k$ where $ r_k= b - Ax_k$ and $e_k = x-x_k$. This leads to the following algorithm.

\begin{alg}\label{alg:ROAP3}(Restarted orthogonally accumulated projection method-ROAP3) Let $A$ in be an unsymmetric and nonsingular matrix and $b\in R^n$ a non-zero vector. Let $\epsilon $ be a given tolerance. The following procedure produces an approximation to the solution $x$ to system (\ref{def:sys}).
   \begin{quote}
      $err = 1, r = b;  x = 0(\in R^n)$\\
      while err $> \epsilon $
      \begin{quote}
      $t = ||A'r||, v_1 = A'r/t, c_1 = b'r/t, u_1 = v_1$ \\
      $ x_1 = c_1 v_1, \beta_0 = \gamma_0 = 0$, $v_0 = u_0 =0(\in R^n)$\\
      for $k =1 $ to $n - 1$\\[-4mm]
        \begin{quote}
           $\alpha_k = u_k' A v_k$ \\
           $ p_k = Av_k - \alpha_k u_k - \beta_{k-1} u_{k-1} $\\
           $\gamma_k = ||p_k||$ \\
           $u_{k+1} = p_k / \gamma_k $ \\
           $q_k = A'u_k - \alpha_k v_k - \gamma_{k-1} v_{k-1} $\\
           $\beta_k = || q_k || $\\
           $ v_{k+1} = q_k / \beta_k $\\
           $c_{k+1} = \beta_k^{-1}(b'u_k - \alpha c_k - \gamma_{k-1} c_{k-1}) $\\
           $\theta = \cos^{-1}(xk'v_k/||x_k||)$\\
           if $| \pi/2-\theta | = 0$\\
               \mbox{~~~~~~} $x_{k+1} = x_k + c_{k+1} v_{k+1}$ \\
           else\\
              \mbox{~~~~~~}$ r = b - Ax_k$\\
              \mbox{~~~~~} break;\\
           end\\
        \end{quote}\vskip -4mm
      end \\ 
      $x = x + x_k$, $r = b- Ax.$\\
       err = $||b - Ax|/||b||$\\
      \end{quote}
      \vskip -4mm
      end     
   \end{quote}
\end{alg}

Remark: It is easy to see that the above restarted orthogonally accumulated projection method is a convergent
iterative scheme since the resulted error vector sequence $e_k$ produced in every restart iteration is a strictly decreasing sequence in terms of their Eucleadean norms.

\section{An orthogonally accumulated projection through bidiagonalization}\label{sec:bidiag}

In this section we propose an iterative scheme similar to the OAP algorithm introduced in section \ref{sec:tridiag}. Instead of using Lanczos-like process based on tridiagonalization of an unsymmetric matrix, we show in this section that an analogous Lanczos-like process can also be based on bidiagonalization of unsymmetric matrix.

\subsection{Matrix bidiagonalization}

Any matrix $A \in R^{n\times m}$ can be transformed into the following bidiagonal form
  \begin{equation} \label{eq:bidiagonal_form1}
     U'AV = T
  \end{equation}
where $T$ is tridiagonal
 $$ T= \left(\begin{array}{cccccc}
            \alpha_1 & \beta_1  & 0       & \cdots&  0 & 0 \\
            0 & \alpha_2 & \beta_2 & \cdots& 0 & 0 \\
            \vdots   & \vdots &  \vdots &   \vdots & \vdots & \vdots  \\
            0        &   0    &  0      &    \cdots & \alpha_{n-1}& \beta_{n-1} \\
             0        &   0    &  0      &    \cdots & 0 & \alpha_n \\
             \end{array} \right),
$$ both $U$ and $V$ are orthogonal, i.e., $U'U = V'V = I_n$. Of course this transform can be accomplished stably  by applying Householder transformations on both sides of $A$.
However a more efficient Lanczos-like process can be used to do the bidiagonalization. To illustrate this we rewrite equation (\ref{eq:bidiagonal_form1}) into the following forms
 \begin{equation} \label{eq:bidiagonal_form2_1}
    A V = U T
   \end{equation}
and
 \begin{equation} \label{eq:bidiagonal_form2_2}
    A' U= V T'
   \end{equation}
Equating $k$-th column of both sides of (\ref{eq:bidiagonal_form2_1}) and (\ref{eq:bidiagonal_form2_2})
we have
 \begin{equation} \label{eq:bidiagonal_form2_3}
    A v_k = \alpha_k u_k +  \beta_{k-1} u_{k-1},  \quad k=1,2,\cdots, n
   \end{equation}
and
 \begin{equation} \label{eq:bidiagonal_form2_4}
    A' u_k = \beta_k v_{k+1} + \alpha_k v_k,  \quad k=1,2,\cdots, n-1
   \end{equation}
with $ \beta_0 =0$, where $v_k$ and $u_k$ denote the $k$-th columns of
matrix $V$ and $U$ separately, $ v_0$ is a zero vector, i.e, $V = [v_1, v_2, \cdots, v_n]$ and
$U = [u_1, u_2, \cdots , u_n]$. Especially we have
\begin{equation} \label{eq:1st_col}
   A v_1 = \alpha_1 u_1  \mbox{~~~ and~~~~ } A' u_1 = \alpha_1 v_1 + \beta_1 v_2
\end{equation}
which suggests that if $v_1$ is given, then $u_1$ and $v_2$ can be calculated successively.
The rest vectors $v_k$($k>3$) and $u_k$ ($k\ge 2$) can be calculated by rewriting (\ref{eq:bidiagonal_form2_3})
and (\ref{eq:bidiagonal_form2_4}) as follows
\begin{equation} \label{eq:bidiagonal_form2_5}
    u_k = \frac{1}{\alpha_k} (A v_k - \beta_{k-1} u_{k-1} )
   \end{equation}
and
 \begin{equation} \label{eq:bidiagonal_form2_6}
   v_{k+1} = \frac{1}{\beta_k} (A' u_k - \alpha_k v_k )
   \end{equation}
The following algorithm depicts the above process.

\begin{alg}\label{alg:OAP2}
   Let $A \in R^{n \times n}$, $v_1 \in R^n$ be a unit vector.
  \begin{quote}
     $\beta_0=0, u_0=0\in R^n$\\
      for $k=1$ to $n-1$ \\
       \mbox{~~~}$w_k = Av_k - \beta_{k-1} u_{k-1}$\\
       \mbox{~~~}$\alpha_k = || w_k||$\\
       \mbox{~~~}$u_k = w_k / \gamma_k$\\
       \mbox{~~~}$q_k = A' u_k - \alpha_k v_k$\\
       \mbox{~~~}$\beta_k = ||q_k||$\\
       \mbox{~~~}$ v_{k+1} = q_k/ \beta_k$\\
      end
  \end{quote}
\end{alg}

In exact arithmetic operations the above Lanczos-like iteration will produce two orthonormal vector sequences $v_1, v_2, \cdots v_n$ and $u_1, u_2, \cdots u_n$ with any starting unit vector $v_1$, assuming no break-down happens(i.e.,  $ \beta_k \neq 0$ for all $k$). Note that in each loop in the iteration one needs only two matrix-vector multiplications as its major flop counts, this makes it very effective when dealing with bidiagonalizations of large and sparse matrices.

\subsection{Orthogonally accumulated projection}

To develop a corresponding accumulated projection method, we need a sequence of orthonormal vectors $\{v_k\}_1^n$ and the inner-product between each $v_k$ and $x$, the exact solution to the system (\ref{def:sys}). Again this can be easily obtained if we choose a starting unit vector $v_1$ with $c_1=x'v_1$
given, since we have by multiplying both sides of equation (\ref{eq:bidiagonal_form2_6}) by $x$
\begin{equation} \label{comput_ci}
    c_{k+1} = x'v_{k+1} = \frac{1}{\beta_k} ( x'A'u_k  - \alpha_k x'v_k) =
    \frac{1}{\beta_k} (b'u_k - \alpha_k c_k ),
\end{equation}
since $Ax = b$.
This implies that if $c_1$ is known, then all the other subsequent $c_i(i=2,3,\cdots n)$ can be calculated by (\ref{comput_ci}).
These process can be described in the following algorithm, which can be viewed as an augumented Lanzcos iteration for solving linear system of equations.

\begin{alg}\label{alg:OAP2}(orthogonally accumulated projection method-OAP2) Let $A$ in be an unsymmetric and nonsingular matrix and $b\in R^n$ a non-zero vector. Let $v_1$ be a unit vector and $c_1( = x'v_1)$ given, where $x$ is the solution to (\ref{def:sys}).
The following process gives the exact solution $x$ to system $Ax = b.$
   \begin{quote}
      $ x_1 = c_1 v_1$,$\beta_0=0$, $ u_0 =0(\in R^n)$\\
      for $k =1 $ to $n - 1$\\[-4mm]
        \begin{quote}

           $ p_k = Av_k  - \beta_{k-1} u_{k-1} $\\
           $\alpha_k = ||p_k||$ \\
           $u_k = p_k / \alpha_k $ \\
           $q_k = A'u_k - \alpha_k v_k  $\\
           $\beta_k = || q_k || $\\
           $ v_{k+1} = q_k / \beta_k $\\
           $c_{k+1} =(b'u_k - \alpha c_k )/\beta_k $\\
           $x_{k+1} = x_k + c_{k+1} v_{k+1}$ \\
        \end{quote}\vskip -4mm
      end
   \end{quote}
\end{alg}
Note that there are only two matrix-vector multiplications involved, and storage for extra three vectors
is needed besides that for the coefficient matrix $A$. Also the flop counts for each oap loop is O(mn) in case $A$ is sparse(having an average of $m$ none-zero elements in each row) and large.

It is also easy to verify the orthonormality of the vector sequences $\{v_i\}_1^n$ and $\{u_i\}_1^n$ in Algorithm \ref{alg:OAP2}, the conclusion is stated in the following.

\begin{thm}
Let $A$ be unsymmetric and nonsingular, $b\in R^n$ and $x$ is the solution to $Ax=b$. The vector sequences $v_k \, (k=1,2,\cdots, n)$ and $u_k\, (k=1,2,\cdots, n)$produced in Algorithm \ref{alg:OAP2} are orthonormal, assuming no breakdown happens, i.e., $\beta_k \neq 0$  for any $k=1,2,3,\cdots, n-1$.
\end{thm}

Proof. Apparently all vectors $v_i$ and $u_i$ ($i = 1, 2, \cdots, n$) are unit vectors.
We first  show that $v_2'v_1=0$ and $u_2 'u_1=0$.

Note that$$
\begin{array}{ll}
     v_2'v_1 = 0 & \Leftrightarrow (A'u_1 - \alpha_1 v_1)' v_1 = 0 \\
                 & \Leftrightarrow \alpha_1 = v_1'A'u_1   \\
                 & \Leftrightarrow \alpha_1 = \alpha_1 u_1 u_1\\
\end{array}
$$
the last equation holds true since $u_1$ is a unit vector. Similarly we have
$$
\begin{array}{ll}
     u_2'u_1 = 0  & \Leftrightarrow (Av_2 - \beta_1 u_1)' u_1 = 0 \\
                  &\Leftrightarrow \beta_1 = v_2'A'u_1 \\
                  &\Leftrightarrow \beta_1=v_2'(\alpha_1 v_1 + \beta_1 v_2)\\
\end{array}
$$
The last equation is true since $v_1'v_2 = 0 $ and $v_2 $ is a unit vector.
By induction, we assume $v_1, v_2, \cdots, v_k$ and $u_1, u_2, \cdots, u_k$ are orthonormal sequences of vectors, we need to show that $v_{k+1}'v_i =0$ and $u_{k+1}'u_i =0$ for $i \le k$.

In fact
$$ \begin{array}{ll}
u_{k+1}' u_k = 0 & \Leftrightarrow (Av_{k+1} -  \beta_k u_k)' u_k =0 \\
           &\Leftrightarrow \beta_k = u_k' A v_{k+1} \\
            &\Leftrightarrow \beta_k = v_{k+1}' (A'u_k)' \\
            &\Leftrightarrow \beta_k = v_{k+1}'( \beta_k v_{k+1} + \alpha_k v_k)
\end{array}$$
        and

For $i < k $ we have
$$
\begin{array}{ll}
    u_{k+1}' u_i =0 &\Leftrightarrow (A v_{k+1} - \beta_k u_k)' u_i =0  \\
                   &\Leftrightarrow  v_{k+1}' (A'u_i) =0 \\
                   &\Leftrightarrow  v_{k+1}' ( \beta_i v_{i+1} + \alpha_i v_i) =0
\end{array}
$$ The last equation holds true since by assumption we have $v_k$ are orthogonal to $v_i$ for any $i<k$.
Similarly one can prove $u_{k+1}'u_i=0$ for $i \le k$.  $\Box$

To handle the issue of loss of orthogonality, a restarted orthogonally accumulated projection can be used, which is analogous to Algorithm \ref{alg:ROAP3} and is stated as

\begin{alg}\label{alg:ROAP2}(Restarted orthogonally accumulated projection method-ROAP2) Let $A\in R^{n\times n}$ in be an unsymmetric and nonsingular matrix and $b\in R^n$ a non-zero vector. Let $\epsilon (<<1) $ be a given tolerance. The following procedure produces an approximation to the solution $x$ to system (\ref{def:sys}).
   \begin{quote}
      $err = 1, r = b;  x = 0(\in R^n)$\\
      while err $> \epsilon $
      \begin{quote}
      $t = ||A'r||, v_1 = A'r/t, c_1 = b'r/t$ \\
      $ x_1 = c_1 v_1, \beta_0 =0$, $u_0 =0(\in R^n)$\\
      for $k =1 $ to $n - 1$\\[-4mm]
        \begin{quote}
           $ p_k = Av_k - \beta_{k-1} u_{k-1} $\\
           $\alpha_k = ||p_k||$ \\
           $u_{k+1} = p_k / \alpha_k $ \\
           $q_k = A'u_k - \alpha_k v_k $\\
           $\beta_k = || q_k || $\\
           $ v_{k+1} = q_k / \beta_k $\\
           $c_{k+1} =(b'u_k - \alpha_k c_k)/\beta_k $\\
           $\theta = \cos^{-1}(xk'v_k/||x_k||)$\\
           if $| \pi/2-\theta | = 0$\\
               \mbox{~~~~~~} $x_{k+1} = x_k + c_{k+1} v_{k+1}$ \\
           else\\
               \mbox{~~~~~} break;\\
           end\\
        \end{quote}\vskip -4mm
      end \\ 
      $x = x + x_k$, $r = b- Ax.$\\
       err = $||b - Ax|/||b||$\\
      \end{quote}
      \vskip -4mm
      end     
   \end{quote}
\end{alg}

\section{Numerical Experiments}

In this section we will examine the numerical behavior of the orthogonally accumulated projection methods proposed in previous sections. OAP methods are used to solve linear system of equations with unsymmetric as well as symmetric coefficient matrices, the results are compared with those obtained by using some benchmark Krylov subspace methods packaged in Matlab. In all the experiments we use the relative residual norm ($||b-Ax_k||/||b||$) as the index for convergence, and the convergence tolerance is set as $10^{-6}$. Also the parameter ``restart" of GMRES is always set as $5$ and parameter ``maximum iteration number" for GMRES is set as the
size of each system in all the experiments.

\textbf{Example 1}. Consider the following convection diffusion problem
$$ \triangle u + p_1 u_x + p_2 u_y + p_3 u = f(x,y)$$
defined on unit square $[0,1]^2$, which usually describes physical phenomena where particles, energy, or other physical quantities are transferred inside a physical system due to two processes: diffusion and convection.  We use the five point finite difference method to discretize the problem, which leads to the following discretized equation
$$ \begin{array}{rl} \frac{2u_{i,j} - u_{i-1, j} - u_{i+1,j} }{(h_x)^2}
      + \frac{2 u_{i,j} - u_{i,j-1} - u_{i,j+1} }{(h_y)^2}
      + p_1 \frac{u_{i+1,j } - u_{i-1,j}}{2h_x} +p_2 \frac{u_{i,j+1} - u_{i,j-1}} {2h_y} + p_3u_{i,j}\\[4mm]
      = f(x_i,x_j)
      \end{array}$$
on each node point $(x_i,y_j)$, where
$u_{i,j} \equiv u(x_i,y_j)$, $h_x$, $h_y$ denote the step size on $x$-axis and $y$-axis direction respectively.  This leads to a linear system of equation $Ax = b$ with $A$ a block tridiagonal unsymmetric matrix.

Table \ref{table:1} shows the comparison of iterative errors among ROAP2 and ROAP3 and some other prevalent Krylov subspace methods. It seems that OAP methods produces better precision than other methods in these experiments, especially than that of GMRES.

\begin{table}[!h]
\tabcolsep 0pt
\caption{Example 1: Comparison of relative errors }
\label{table:1}
\vspace*{-12pt}
\begin{center}
\def\temptablewidth{0.9\textwidth}
{\rule{\temptablewidth}{1pt}}
\begin{tabular*}{\temptablewidth}{@{\extracolsep{\fill}}|c|c|c|c|c|c|c|}
    n & ROAP2 & ROAP3  & GMRES & LSQR & QMR & BiCG  \\   \hline
  90   & 6.0659e-12& 4.8411e-8 &5.8966e-7 &1.1206e-7 &1.4894e-7& 5.9289e-8\\ \hline
  171  & 6.1516e-9& 8.8727e-8 &1.0452e-6& 1.2838e-7 &1.5546e-7 &5.8738e-8\\ \hline
  361  & 7.3004e-8& 1.8632e-8 &6.5894e-7& 4.1377e-8& 3.9033e-8 &3.0593e-8\\ \hline
  551  & 1.2491e-10& 1.5095e-8 &1.0729e-6& 5.3868e-8& 1.1865e-7 &6.4840e-8\\ \hline
  741  & 9.0775e-10& 4.3456e-9& 1.1596e-6& 5.8852e-8& 1.1784e-7& 4.2417e-8\\ \hline
  1131 & 2.7517e-9& 1.9215e-8& 1.1474e-6& 6.1654e-8 &9.0458e-8& 3.4545e-8\\ \hline
  1521 & 1.3374e-8& 2.4574e-8& 1.1846e-6 &2.2966e-8 &5.6823e-8& 1.9139e-8 \\ \hline
  2401 & 5.0582e-9& 7.3975e-9 &1.2118e-6& 2.4055e-8& 5.0915e-8 &2.0182e-8 \\ \hline
       \end{tabular*}
       {\rule{\temptablewidth}{1pt}}
       \end{center}
       \end{table}

\begin{table}[!h]
\tabcolsep 0pt
\caption{Example 1: Comparison of iteration numbers }
\label{table:2}
\vspace*{-12pt}
\begin{center}
\def\temptablewidth{0.8\textwidth}
{\rule{\temptablewidth}{1pt}}
\begin{tabular*}{\temptablewidth}{@{\extracolsep{\fill}}|c|c|c|c|c|c|c|}
  n & ROAP2 & ROAP3  & GMRES & LSQR & QMR & BiCG  \\   \hline
  90 &  2 &6& 4 &77 &27& 28\\   \hline
  171&  2& 6& 9& 178& 44& 46\\   \hline
  361& 5 &6 &12 &188& 47& 47\\   \hline
  551 &3& 12& 19& 479& 68& 70\\   \hline
  741 &2& 10& 26& 744 &86& 90\\   \hline
  1131& 6 &8& 35& 917 &93& 96\\   \hline
  1521& 9& 8& 43& 764& 94 &96\\   \hline
  2401& 1 &8& 66 &1190& 118 &120\\   \hline
\end{tabular*}
{\rule{\temptablewidth}{1pt}}
\end{center}
\end{table}

\begin{figure}[htp]
   \subfigure[distribution of non zero elements of A ]{\includegraphics[width=6cm]{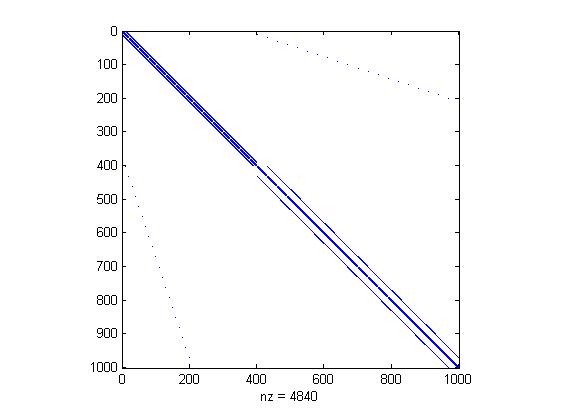}
   \label{fig:c} }
  \hspace{5mm}
   \subfigure[distribution of non zero elements of A ]{\includegraphics[width=6cm]{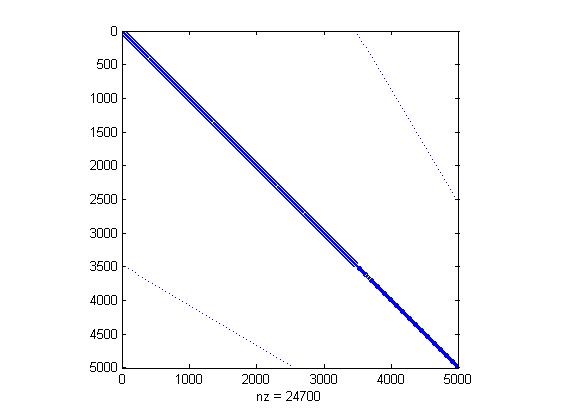}
   \label{fig:d} }
   \caption{Example 2: Pattern of non-zero elements distribution  }
\end{figure}

\textbf{Example 2.} We test the Poisson problem defined on a L-shaped domain  $ [0,1]\times[0,\frac{1}{2}] \cup [0,\frac{1}{2}]  \times[\frac{1}{2},1]$. The resulted  coefficient matrices are symmetric and positive definit. They ususlly have zero pattern shown as in Figure \ref{fig:c} and Figure \ref{fig:d}. The comparison of relative errors among OAP and other Krylov subspace methods are shown in Table \ref{table:3}. It seems that again OAP methods produce better precision than other methods in terms of relative errors.

\begin{table}[!h]
\tabcolsep 0pt
\caption{Example 2: Comparison of relative errors }
\label{table:3}
\vspace*{-12pt}
\begin{center}
\def\temptablewidth{\textwidth}
{\rule{\temptablewidth}{1pt}}
\begin{tabular*}{\temptablewidth}{@{\extracolsep{\fill}}|c|c|c|c|c|c|c|c|c|c|}
 n& ROAP2&  ROAP3& PCG& GMRES&  LSQR&  QMR&  BiCG&  SYMMLQ&  MINRES  \\ \hline
 200& 1.1714e-7 & 9.9202e-9&  7.3111e-8&  7.1231e-7&  7.7353e-8 & 1.1013e-7&  7.3111e-8  & 7.3111e-8  & 1.1013e-7\\ \hline
 500& 1.5743e-7 &  1.072e-7 & 2.3488e-7 &  4.608e-6 & 4.2399e-7&  7.5256e-7 & 2.3488e-7  & 2.3488e-7&
  7.5256e-7\\ \hline
 1000& 2.9599e-7&  1.1256e-7 & 4.0789e-7&  8.6828e-6  &7.2121e-7 & 1.2611e-6 & 4.0789e-7  & 4.0789e-7&
  1.2611e-6\\ \hline
 1400&  3.4842e-7 &  2.274e-7 & 6.8198e-7&  1.1024e-5&  4.9706e-7 & 1.8123e-6&  6.8198e-7 & 6.8198e-7&
  1.8123e-6\\ \hline
1700&   2.4713e-7 & 5.1705e-7 &  5.885e-7&  1.3928e-5  &7.4723e-7&  3.6628e-6&   5.885e-7   & 5.885e-7&
  3.6628e-6\\ \hline
 2100 &2.5727e-8 & 7.8073e-7 & 5.6357e-7  &1.7001e-5 & 2.8269e-7 & 1.0733e-6 & 5.6357e-7  & 5.6357e-7&
  1.0733e-6 \\ \hline
\end{tabular*}
{\rule{\temptablewidth}{1pt}}
\end{center}
\end{table}

\begin{table}[!h]
\tabcolsep 0pt
\caption{Example 2: Comparison of iteration numbers }
\label{table:4}
\vspace*{-12pt}
\begin{center}
\def\temptablewidth{.7\textwidth}
{\rule{\temptablewidth}{1pt}}
\begin{tabular*}{\temptablewidth}{@{\extracolsep{\fill}}|c|c|c|c|c|c|c|c|c|c|}
 n& ROAP2&  ROAP3& PCG& GMRES&  LSQR&  QMR&  BiCG&  SYMMLQ&  MINRES  \\ \hline
 200& 6 &6 &41& 8& 170& 41& 41 &40 &41\\ \hline
 500&6 &12 &69 &14& 425& 68& 69& 68& 68\\ \hline
 1000& 13& 27& 94& 23& 826& 92& 94 &93& 92\\   \hline
 1400&9& 42& 111& 34& 1149 &108& 111& 110& 108\\ \hline
 1700& 7 &56& 120 &39& 1381& 115 &120 &119& 115\\ \hline
 2100&6& 62& 111 &44& 1052& 109 &111 &110& 109\\ \hline

\end{tabular*}
{\rule{\temptablewidth}{1pt}}
\end{center}
\end{table}

\textbf{Example 3} We take unsymmetric tridiagonal matrix $A=diag\{-1,2,-1.1\}_n$ as coefficient matrix, and the right hand vector $b$ is taken such that the exact solution is a vector contains the function values of $x(t)=t(1-t)e^t$ at grid points $t=h:h:1-h$, where $h=1/n$.
  The relative errors and iterative numbers resulted from using OAP and other Krylov subspace methods are shown in the Table \ref{table:5} and Table \ref{table:6} respectively. Note that the coefficient matrix has very large condition number as $n$ increases, and the condition numbers are listed in the second column in Table \ref{table:6}.
 \begin{table}[!h]
\tabcolsep 0pt
\caption{Example 3: Comparison of relative errors }
\label{table:5}
\vspace*{-12pt}
\begin{center}
\def\temptablewidth{0.9\textwidth}
{\rule{\temptablewidth}{1pt}}
\begin{tabular*}{\temptablewidth}{@{\extracolsep{\fill}}|c|c|c|c|c|c|c|}

  n & ROAP2 & ROAP3  & GMRES & LSQR & QMR & BiCG  \\   \hline
  600 & 3.0413e-4 &3.0413e-4 & 1.0063e-3& 3.0414e-4& 9.8330e-4 & 1.1523e-3\\   \hline
  900 & 1.6567e-4 &1.6567e-4 &5.3850e-4 &1.6569e-4& 5.2552e-4 &6.1508e-4\\   \hline
  1200&  1.0765e-4& 1.0765e-4& 3.4693e-4& 1.0767e-4 &3.3994e-4& 3.9558e-4\\   \hline
  1500 & 7.7045e-5 &7.7045e-5& 2.4720e-4& 7.7080e-5& 2.4247e-4& 2.8136e-4\\   \hline
  1800& 5.8620e-5& 5.8620e-5 &1.8768e-4& 5.8666e-5 &1.8396e-4& 2.1319e-4\\   \hline
  2100 & 4.6524e-5 &4.6524e-5& 1.4885e-4& 4.6581e-5& 1.4553e-4& 1.6869e-4\\   \hline

       \end{tabular*}
       {\rule{\temptablewidth}{1pt}}
       \end{center}
       \end{table}

\begin{table}[!h]
\tabcolsep 0pt
\caption{Example 3: Comparison of iteration numbers }
\label{table:6}
\vspace*{-12pt}
\begin{center}
\def\temptablewidth{0.8\textwidth}
{\rule{\temptablewidth}{1pt}}
\begin{tabular*}{\temptablewidth}{@{\extracolsep{\fill}}|c|c|c|c|c|c|c|c|}
  n & cond(A)& ROAP2 & ROAP3  & GMRES & LSQR & QMR & BiCGstab  \\   \hline
  600 & 3.8846e+14    & 6 &6 &600 &428 &302& 23\\   \hline
  900  & 1.2466e+21   & 6 &6& 900& 388 &370& 23\\   \hline
  1200 &  3.6164e+27  & 6 &6 &1200 &357 &354 &23\\   \hline
  1500 & 1.8172e+33   &  6& 5 &1800 &316& 327& 23\\   \hline
  1800& 2.6357e+39   & 6 &5& 1800 &316& 327& 23\\   \hline
  2100 & 8.0531e+45&  6& 5 &2100& 296& 311& 23 \\   \hline
\end{tabular*}
{\rule{\temptablewidth}{3pt}}
\end{center}
\end{table}

\textbf{Example 4} We use Matlab  routine $rand()$ to produce coefficient matrix  $A$, the right hand side vector $b$ is taken so that the exact solution is a vector contains the function values of $x(t)=t(1-t)e^{3t}$ at grid points $t = i*h\,(i=1,2,\cdots, n)$, where $h=1/n$.
  The relative errors and iterative numbers resulted from using OAP and other Krylov subspace methods are shown in the Table \ref{table:7} and Table \ref{table:8} respectively.  We found that except LSQR, other tested methods such as QMR,BiCG, BiCGstab and GMRES all fail to produce convergent resultus in these experiments.

 \begin{table}[!h]
\tabcolsep 0pt
\caption{Example 4: Comparison of relative residual }
\label{table:7}
\vspace*{-12pt}
\begin{center}
\def\temptablewidth{0.9\textwidth}
{\rule{\temptablewidth}{1pt}}
\begin{tabular*}{\temptablewidth}{@{\extracolsep{\fill}}|c|c|c|c|c|c|c|}
    n & ROAP2 & ROAP3  & GMRES & LSQR & QMR & BiCGstab  \\   \hline
  300& 9.9465e-7 &7.5874e-7 &2.0905e-2 &9.9106e-7 &7.6749e-3 &2.1436e-2\\   \hline
  600& 7.6515e-7 &5.7610e-7& 1.4925e-2& 9.9968e-7 &1.5316e-2 &1.5317e-2\\   \hline
  900& 5.0974e-7 &8.0796e-7 &1.1373e-2 &9.9755e-7 &1.1457e-2 &1.1468e-2\\   \hline
       \end{tabular*}
       {\rule{\temptablewidth}{1pt}}
       \end{center}
       \end{table}

\begin{table}[!h]
\tabcolsep 0pt
\caption{Example 4: Comparison of iteration numbers }
\label{table:8}
\vspace*{-12pt}
\begin{center}
\def\temptablewidth{0.8\textwidth}
{\rule{\temptablewidth}{1pt}}
\begin{tabular*}{\temptablewidth}{@{\extracolsep{\fill}}|c|c|c|c|c|c|c|}
  n & ROAP2 & ROAP3  & GMRES & LSQR & QMR & BiCG  \\   \hline
  300 & 106& 229 &300& 536& 1498& 1\\   \hline
  600 & 15& 52 &600& 985 &1 &1 \\   \hline
  900 & 20& 30 &900& 866& 2& 1\\   \hline
\end{tabular*}
{\rule{\temptablewidth}{1pt}}
\end{center}
\end{table}

\section{Comments and Summary}

The OAP methods introduced in this paper still belong to the category of extended Krylov subspace methods since they rely on the construction of Krylov subspaces $K_m(A,v)$ and $K_m(A',v)$ with fixed coefficient matrix. Although they are also derived from Lanczos process, just like some other Krylov subspace methods such as QMR, BiCG, BiCGstab, MINRES, CG; a major feature that makes OAP different than the other methods is the detection of   loss of orthogonality is used in OAP, while the others usually do nothing to deal with loss of orthogonality. This might be the explanation of the instability of these classical Krylov subspace methods. Also it is easy to show the restart strategy used in OAP leads to a convergent iterative scheme, while restarted GMRES does not always guarantee a convergent process. As a matter of fact, it can be shown that CG can be viewed as a generalized OAP method where the orthogonality between vectors $v_1$ and $v_2$ is defined as $v_1'Av_2=0$ instead of $v_1'v_2=0$, thus a restart CG method can also be derived and is also convergent, while successful adoptionof restart strategy( which leads to a convergent iterative scheme) on other classical Krylov subspace methods are hard.

It is also possible for us to develop accelerative schemes similar with those presented in\cite{pengPAP} \cite{pengSAP} for OAP algorithms in case of very large scale computation.

\bibliographystyle{plain}
\bibliography{peng_refv2}

\end{document}